%% file: yoneda.tex
\documentclass[12pt]{amsart}
\usepackage{a4wide}
\usepackage{eucal}
\usepackage{umlaut}
\usepackage{latexsym}
\usepackage{amssymb}
\usepackage{verbatim}
\usepackage[active]{srcltx}
\usepackage[all]{xy}
\usepackage{listings}
\usepackage{rotating}
\usepackage{tikz}
\usepackage{multirow}
\usepackage[pdfauthor   = {Mohamed\ Barakat\ and\ Barbara\ Bremer},
            pdftitle    = {Higher\ Extension\ Modules\ and\ the\ Yoneda\ Product},
            pdfsubject  = {13D02; 13D05; 13D07; 16E05; 16E10; 16E30},
            pdfkeywords = {Yoneda\ products;\ higher\ extension\ modules},
            bookmarks=false,
            colorlinks=true,
            pagebackref=true,
            hyperindex=true,
            linkcolor=blue,
            citecolor=blue,
            urlcolor=blue,
            ps2pdf=true
            ]{hyperref}

\input amathmoe.sty

\input pre.tex

\begin{document}

\lstset{basicstyle=\tiny,
        frame=single,
        stringstyle=\ttfamily,
        showstringspaces=true}

\author{Mohamed Barakat \& Barbara Bremer}
\address{Lehrstuhl B f\"ur Mathematik, RWTH-Aachen University, 52062 Germany}
\email{\href{mailto:Mohamed Barakat <mohamed.barakat@rwth-aachen.de>}{mohamed.barakat@rwth-aachen.de}, \href{mailto:Barbara Bremer <barbara.bremer@rwth-aachen.de>}{barbara.bremer@rwth-aachen.de}}

\PUSH{Yoneda.tex}%
\input Yoneda.tex%
\POP

\PUSH{yoneda.bbl}%
\input yoneda.bbl
\POP




\end{document}

%% file: pre.tex
\DeclareMathOperator{\Hom}{Hom}
\DeclareMathOperator{\Ext}{Ext}
\DeclareMathOperator{\Yon}{Yon}
\DeclareMathOperator{\YExt}{YExt}
\DeclareMathOperator{\ExtMod}{ExtMod}

\DeclareMathOperator{\coker}{coker}
\DeclareMathOperator{\img}{im}
\DeclareMathOperator{\defect}{def}

\DeclareMathOperator{\RR}{R}
\DeclareMathOperator{\Sat}{S}

\newcommand\id{\mathrm{id}}
\renewcommand\phi{\varphi}

\newcommand\homalg{\mathtt{homalg}}

\newcommand\tE{\mathtt{E}}

\newcommand\tH{\mathtt{H}}
\newcommand\tK{\mathtt{K}}
\newcommand\tL{\mathtt{L}}
\newcommand\tM{\mathtt{M}}
\newcommand\tN{\mathtt{N}}

\newcommand\cocoa{{\hbox{\rm C\kern-.13em o\kern-.07em C\kern-.13em o\kern-.15em A}}}

\newcommand{\lklammer}[1]{      
  \ensuremath{\left\{ \rule{0pt}{#1} \right.}
}
\newcommand{\rklammer}[1]{
  \ensuremath{\left\} \rule{0pt}{#1} \right.}
}

\newcommand{\jumpdir}[2]{\save[]+#2*{#1} \restore}

\newcommand{\encircle}[3]{
  \setlength{\unitlength}{1cm}
  \begin{picture}(0,0)
    \put(#1,#2){\circle{#3}}
  \end{picture}
}

\newcommand{\enbox}[4]{
  \setlength{\unitlength}{1cm}
  \begin{picture}(0,0)
    \put(#1,#2){\framebox(#3,#4){}}
  \end{picture}
}

\newcommand{\enboxdash}[4]{
  \setlength{\unitlength}{1cm}
  \begin{picture}(0,0)
    \put(#1,#2){\dashbox{0.1}(#3,#4){}}
  \end{picture}
}
\newcommand{\kommentar}[1]{}

%% file: Yoneda.tex
\title{Higher Extension Modules and the {\sc Yoneda} Product}
\date{2008}

\begin{abstract}
A chain of $c$ submodules $E =: E_0 \geq E_1 \geq \cdots \geq E_c \geq E_{c+1} := 0$ gives rise to $c$ composable $1$-cocycles in $\Ext^1(E_{i-1}/E_i,E_i/E_{i+1})$, $i=1,\ldots,c$. In this paper we follow the converse question: When are $c$ composable $1$-cocycles induced by a module $E$ together with a chain of submodules as above? We call such modules $c$-extension modules. The case $c=1$ is the classical correspondence between $1$-extensions and $1$-cocycles. For $c=2$ we prove an existence theorem stating that a $2$-extension module exists for two composable $1$-cocycles $\eta^M_L\in\Ext^1(M,L)$ and $\eta^L_N\in\Ext^1(L,N)$, if and only if their {\sc Yoneda} product $\eta^M_L\circ\eta^L_N\in\Ext^2(M,N)$ vanishes. We further prove a modelling theorem for $c=2$: In case the set of all such $2$-extension modules is non-empty it is an affine space modelled over the abelian group that we call the first extension group of $1$-cocycles, $\Ext^1(\eta^M_L,\eta^L_N) :=\Ext^1(M,N)/(\Ext^0(M,L)\circ \eta^L_N+\eta^M_L\circ\Ext^0(L,N))$.
\end{abstract}

\maketitle

\section{Introduction}

Let $D$ be a ring with one. We will not restrict this condition further in the paper, only for the examples we need $D$ to be a ring over which one can effectively compute extension groups.

We are interested in the computational problem of (re)constructing a (left) $D$-module from given subfactor modules: \\
A chain of $c$ submodules $E =: E_0 \geq E_1 \geq \cdots \geq E_c \geq E_{c+1} := 0$ gives rise to $c+1$ subfactor modules $E_i/E_{i+1}$, for $i=0,\ldots,c$ and we want to describe the additional data necessary to reconstruct $E$ out of these $c+1$ subfactors.

Since such a chain describes a multiple extension process, we call the module $E$ together with its chain of submodules a \emph{$c$-extension module}. This definition appears in Section~\ref{extmodules} together with the adequate notion of equivalence.

It is obvious that the data needed to reconstruct a chain of length $c$ contains the corresponding data of all its subchains. For example, each pair of successive subfactors leads to the subchain $E_i/E_{i+2} \geq E_{i+1}/E_{i+2} \geq 0$, which corresponds to the short exact sequence $0 \leftarrow E_i/E_{i+1} \leftarrow E_i/E_{i+2} \leftarrow E_{i+1}/E_{i+2} \leftarrow 0$. Thus, one is lead to consider the problem for $c=1$ first. So let $E \geq N \geq 0$ be a $1$-extension module with corresponding short exact sequence $0\leftarrow M \leftarrow E \leftarrow N \leftarrow 0$ for $M:=E/N$. Such short exact sequences with fixed factor module $M$ and submodule $N$ are classified by the first extension group\footnote{For simplicity we will ommit the subscript $D$ in $\Ext_D^c(M,N)$.} $\Ext^1(M,N)$. Its elements are called \emph{$1$-extension cocycles}. In particular, for the case $c=1$, the $1$-cocycles in $\Ext^1(M,N)$ are precisely the data describing how to put $M$ on top of $N$ to reconstruct $E$. This is part of the classical {\sc Yoneda} equivalence which gives an alternative way of describing the extension groups $\Ext^c(M,N)$ as the set of all exact sequences $0 \leftarrow M \leftarrow G_0 \leftarrow \cdots \leftarrow G_{c-1} \leftarrow N \leftarrow 0$ up to an appropriate equivalence relation (cf.~Appendix~\ref{composition}). Turning back to our original problem one observes that any $c$-extension module induces such an exact sequence with $G_i:=E_i/E_{i+2}$. However, unlike the case $c=1$, $E_0=E$ is not one of the modules in the exact sequence and the {\sc Yoneda} equivalence cannot be exploited in the same way for $c>1$ as for $c=1$.

We now set $L_i:=E_i/E_{i+1}$ to simplify the notation and proceed inductively by making use of the data describing the subchains $0 \leftarrow L_i \leftarrow G_i \leftarrow L_{i+1} \leftarrow 0$ of length $1$, namely their $1$-cocycles $\eta^{L_i}_{L_{i+1}}\in\Ext^1(L_i,L_{i+1})$. In other words, we want to classify all $c$-extension modules $E$ with $c+1$ prescribed subfactors $L_i$ together with a $1$-cocycle $\eta^{L_i}_{L_{i+1}}$ for each consecutive pair $(L_i,L_{i+1})$. One first observes that not every choice of such $1$-cocycles leads to a $c$-extension module, so it is natural to ask about the necessary and sufficient conditions for such a module with a prescribed $c$-tuple ($\eta^{L_0}_{L_1},\ldots,\eta^{L_c}_{L_{c+1}}$) to exist.

To describe a necessary condition we recall the fact that two compatible short exact sequences $\eta^{L_i}_{L_{i+1}}:0 \leftarrow L_i \leftarrow G_i \leftarrow L_{i+1} \leftarrow 0$ and $\eta^{L_{i+1}}_{L_{i+2}}:0 \leftarrow L_{i+1} \leftarrow G_{i+1} \leftarrow L_{i+2} \leftarrow 0$ can be spliced together to a $2$-extension $0 \leftarrow L_i \leftarrow G_i \leftarrow G_{i+1} \leftarrow L_{i+2} \leftarrow 0$, which we will denote by $\eta^{L_i}_{L_{i+1}}\circ\eta^{L_{i+1}}_{L_{i+2}}$. In Section~\ref{case2} we prove an existence theorem stating that a $2$-extension module exists, \emph{if and only if} $\eta^{L_0}_{L_1}\circ\eta^{L_1}_{L_2}$ is trivial as a $2$-extension. Then it is easy to isolate the vanishing of $\eta^{L_i}_{L_{i+1}}\circ\eta^{L_{i+1}}_{L_{i+2}}$ for all consecutive pairs as a necessary condition for a $c$-extension module to exist (cf.~Section~\ref{3 and higher}). Example~\ref{example not sufficient} shows that this condition is not sufficient for $c>2$.

Although the {\sc Yoneda} equivalence cannot be exploited directly for $c>1$ as for $c=1$ as explained above, it still proves crucial: As a natural transformation it also provides a way to reinterpret the morphisms in the long exact $\Ext$-sequences using the so-called \emph{{\sc Yoneda} product}. This is the content of \cite[Theorem~III.9.1]{ML} in {\sc MacLane}'s book, which we will recall in Section~\ref{yoneda}. It will provide the final step for the proof of the existence theorem in Section~\ref{case2}.

In Section~\ref{coord} we will express the intrinsic approach used so far in a well-known setup which makes explicit computations possible. Using this explicit language we were able to prove a modelling theorem for the case $c=2$ which concluded Section~\ref{case2}: In case the set of all $2$-extension modules is non-empty, it is an affine space modelled over the abelian group $\Ext^1(\eta^M_L,\eta^L_N) :=\Ext^1(M,N)/(\Ext^0(M,L)\circ \eta^L_N+\eta^M_L\circ\Ext^0(L,N))$, which we call the \emph{first extension group of $1$-cocycles}.

The explict language of Section~\ref{coord} allows us to identify the set of all $c$-extension modules for prescribed $1$-cocycles ($\eta^{L_0}_{L_1},\ldots,\eta^{L_c}_{L_{c+1}}$) as the solution space of a system of equations over the ring $D$. The central observation of the whole Section is that it is possible to isolate a certain subsystem, which can be solved independently. The solution of this subsystem has an independent meaning as it corresponds to computing a lift of the first cocycle $\eta^{L_0}_{L_1}$. The system is affine for $c=2$ and this observation shows that it is also triangular. For $c=3$ the system is quadratic, but due to the observation it can be reduced to solving two affine systems. This puts us in the position to easily construct the counter-example~\ref{example not sufficient} mentioned above. For $c>3$ the system is still quadratic and we were not able to reduce it further.

Section~\ref{expository} is provided to fix the notation rather than exposing the standard material which can be found in \cite{HS,ML,Wei}, for example, and will be summarized in the Appendix.

\medskip
We are not aware of direct contributions to the existence problem of higher extension modules in this general frame work. Nevertheless, the work \cite{pr} of {\sc van der Put} and {\sc Reversat} on the {\sc Galois} theory of $q$-difference equations addresses the problem from a more algebraic geometric point of view. The connection to their approach will be subject of future joint work with {\sc van der Put}. In computational group theory multiple extensions have been crucial in {\sc Plesken}'s soluble quotient algorithm \cite{plsq}.

We want to emphasize that we were guided by examples computed using $\homalg$ \cite{BR}, which has been extended by the second author in her Diploma Thesis \cite{bredipl} to provide procedures for the computation of {\sc Yoneda} products and the {\sc Yoneda} equivalence between $c$-extensions and $c$-cocycles. Based on the results of this paper $\homalg$ provides procedures to compute $c$-extension modules for $c=2,3$ (over computable commutative rings) if they exist. Examples of such computations are given in Section~\ref{examples}. The details of these computations and the thesis of the second author will appear on the homepage of $\homalg$ \cite{homalg}.

In recent years it became clear (see for example \cite{Ob,Fl,Mounier,Zerz}, to name a few) how local\footnote{As opposed to what we call \emph{global} control theory, where boundary conditions must be considered.} linear control theory can be rephrased in the language of modules over various specific rings, allowing an extensive use of the homological machinery in \cite{Q,PQ99,CQR05,QR08} for example. The results in this paper should be useful to analyse but also to construct control systems with specific properties.

Finally, the reader is encouraged to follow the line of arguments on the simple Example~\ref{example most simple}. It is important to note that we apply morphisms of left modules from the right. This leads to the use of the row convention for matrices.

\section{$c$-Extensions and the {\sc Yoneda} Composite}\label{expository}

For two $D$-modules $M$, $N$ and a natural number $c$, a \emph{$c$-extension} of $N$ with $M$ is an
exact sequence starting at $N$ and ending in $M$ and running through $c$
intermediate modules. 
$$\xymatrix{
  0 & M \ar[l] & G_0 \ar[l]& \ldots \ar[l] & G_{c-1} \ar[l] & N \ar[l] & 0.
\ar[l]
}$$
Motivated by the homomorphism theorem one can illustrate this exact sequence by indicating maps between the submodule lattices of the different $G_i$'s:
\begin{equation}\tag{Staircase}\label{Staircase}
\xymatrix@C=1.2cm@R=0.3cm{
  *=0{} \ar@{-}[d] \jumpdir{M}{/:a(180) +0.2cm/} & *=0{} \ar[l]
  \ar@{-}[d] \jumpdir{G_0}{/:a(180) +0.2cm/} &  &  \\
  *=0{} & *=0{\mbox{\tiny\textbullet}} \ar[l] \ar@{-}[d] \jumpdir{ L_1 \big\{ } {/:a(-45) +0.5cm/} &  \hspace{.8cm}\ar[l] && \\
  & *=0{} & \hspace{0.1cm}\cdots\hspace{0.1cm} \ar[l] & *=0{} \ar[l]\ar@{-}[d] \jumpdir{G_{c-1}}{/:a(180) +0.2cm/} \jumpdir{ \big\} L_{c-1}  } {/:a(65) +0.7cm/}\\
  & & \hspace{.8cm} & *=0{\mbox{\tiny\textbullet}} \ar@{-}[d] \ar[l] & *=0{} \ar[l] \ar@{-}[d]
  \jumpdir{N}{/:a(180) +0.2cm/}\\
  & & & *=0{} & *=0{} \ar[l] \\
}
\end{equation}
We abreviate the sequence by writing $0\leftarrow M \leftarrow (G_i) \leftarrow N \leftarrow 0$. A $1$-extension is nothing but a short exact sequence.

A $c$-extension $G: 0\leftarrow M \leftarrow (G_i) \leftarrow L \leftarrow 0$ and a $c'$-extension $G': 0\leftarrow L \leftarrow (G'_j) \leftarrow N \leftarrow 0$
may be \emph{spliced} together by taking the composite map
$\xymatrix{G_{c-1} & \ar[l] L & \ar[l]  G'_0 \ar@/^1pc/[ll]}$ to give a
$(c+c')$-extension of $M$ with $N$:
$$
 0 \leftarrow M \leftarrow (G_i) \leftarrow (G'_j) \leftarrow N \leftarrow 0
$$
This is called the \emph{{\sc Yoneda} composite} of $G$ and $G'$ and denoted
$G\circ G'$.

The \emph{{\sc Yoneda} equivalence}
\[
  \Yon:\Ext \rightarrow \YExt
\]
establishes the one-to-one correspondence between elements of the $c$-th extension group and $c$-extensions (cf.~Appendices~\ref{composition},~\ref{equivalence}). Elements of $\Ext^c(M,N)$ are called \emph{$c$-extension cocycles} and we call them \emph{$c$-cocycles} for short. Now one can use the {\sc Yoneda} equivalence to carry the {\sc Yoneda} composite over to cocycles. This is called the \emph{{\sc Yoneda} product of cocycles}:
\[
  \eta \circ \eta' := \Yon^{-1}\left(\Yon(\eta)\circ\Yon(\eta')\right),\quad \eta\in\Ext^c(M,L) \mbox{ and } \eta'\in\Ext^{c'}(L,N).
\]
A shorter way to compute the {\sc Yoneda} product can be found in Appendix~\ref{lift}.

We further define $L_i:=\coker(G_i \leftarrow G_{i+1})$ for $i=0,\ldots,c-1$, with $G_c:=N=:L_c$ (see the (\ref{Staircase}) diagram above). Then, the $c$-extension
\[
  \xymatrix@1{
    G: 0 & \ar[l] M & \ar[l] G_0 & \ar[l] \cdots & \ar[l] G_{c-1} & \ar[l] N & \ar[l] 0
  }
\]
is the {\sc Yoneda} composite of the $c$ short exact sequences
$0 \leftarrow L_i \leftarrow G_i \leftarrow L_{i+1} \leftarrow 0$ for $i=0,\ldots,c-1$. We denote the corresponding $1$-cocycles by $\eta^{L_i}_{L_{i+1}}$.

\section{The {\sc Yoneda} Product and the Connecting Homomorphism}\label{yoneda}

Let $\eta^M_N$ be the corresponding cocycle to the short exact sequence
\[
  0 \leftarrow M \xleftarrow{\pi^G} G \xleftarrow{\iota_G} N \leftarrow 0
\]
and $L$ another module. Then the sequences
\[
  \xymatrix@1{
    \cdots \Ext^{c-1}(L,M) \ar[r]^>(.8){\delta^{c-1}} & \Ext^c(L,N) \ar[r]^{{\iota_G}_*}
    & \Ext^c(L,G) \ar[r]^{{\pi^G}_*} & \Ext^c(L,M) \ar[r]^>(.8){\delta^c} &
    \Ext^{c+1}(L,N) \cdots
  }
\]
and
\[
  \xymatrix@1{
    \cdots \Ext^{c-1}(N,L) \ar[r]^>(.8){\Delta^{c-1}} & \Ext^c(M,L) \ar[r]^{{\pi^G}^*} &
    \Ext^c(G,L) \ar[r]^{{\iota_G}^*} &
    \Ext^c(N,L) \ar[r]^>(.8){\Delta^c} & \Ext^{c+1}(M,L) \cdots
  }
\]
are exact. They start at the left with $0 \to \Hom(L,N) = \Ext^0(L,N)$ and
with $0 \to \Hom(M,L)$, respectively, and continue to the right for all $c \geq 0$. The connecting homomomorphisms turn out to \emph{coincide} with the {\sc Yoneda}
products with $\eta^M_N$:
\[
  \delta^c : \Ext^c(L,M) \to \Ext^{c+1}(L,N), \eta^L_M \mapsto \eta^L_M \circ
  \eta^M_N, 
\]
\[
  \Delta^c : \Ext^c(N,L) \to \Ext^{c+1}(M,L), \eta^N_L \mapsto \eta^M_N \circ
  \eta^N_L.
\]
This is the content of \cite[Theorem~III.9.1]{ML}. The idea is simply the universality of both the connecting homomorphism and the {\sc Yoneda} product. Furthermore
\[
\begin{array}{crclcrcl}
  & {\iota_G}_*(-) & = & - \circ \iota_G & \quad &
  {\pi^G}_*(-) & = & - \circ \pi_G \\
  \mbox{and} \quad & {\iota_G}^*(-) & = & \iota_G \circ - & &
  {\pi^G}^*(-) & = & \pi_G \circ -
\end{array}
\]

For a $k$-cocycle $\eta^M_N\in\Ext^k(M,N)$ one can also define the so called \emph{iterated connecting homomorphisms} using the {\sc Yoneda} product:
\[
  \delta^c_k : \Ext^c(L,M) \to \Ext^{c+k}(L,N), \eta^L_M \mapsto \eta^L_M \circ
  \eta^M_N, 
\]
\[
  \Delta^c_k : \Ext^c(N,L) \to \Ext^{c+k}(M,L), \eta^N_L \mapsto \eta^M_N \circ
  \eta^N_L.
\]

\section{$c$-Extension Modules}\label{extmodules}

We call a module $E$ together with a chain of $c$ submodules
$E \geq E_1 \geq \cdots \geq E_c \geq 0$ a \emph{$c$-extension module}\footnote{This is a special case of a filtered module over a ring with the trivial filtration $D=D_0=D_i$.}. 

For a $c$-extension module $E=:E_0 \geq E_1 \geq \cdots \geq E_c \geq E_{c+1}
:= 0$, $c \geq 1$, we denote by $G(E)$ the $c$-extension
\[
  G(E): \xymatrix@1{
  0 & \ar[l] M & \ar[l] G_0 & \ar[l] \cdots & \ar[l] G_{c-1} & \ar[l] N & \ar[l] 0,
  }
\]
where $M:=E_0/E_1$, $G_i:=G_i(E):=E_i/E_{i+2}$ and $N:=E_c$, together with the natural maps $E_i/E_{i+2} \to E_{i-1}/E_{i+1}$. We say that a $c$-extension $G \in\YExt_D^c(M,N)$ is induced by a $c$-extension module, if $G=G(E)$ for some $c$-extension module $E$:

\begin{equation}\tag{$\ExtMod$}
  \xymatrix@C=1.2cm@R=0.3cm{
    *=0{} \ar@{-}[d] \jumpdir{M}{/:a(180) +0.2cm/} & *=0{} \ar[l]
    \ar@{-}[d] \jumpdir{G_0}{/:a(180) +0.2cm/} &  &  & & &
    *=0{\mbox{\tiny\textbullet}} \ar@{-}[d]\\
    *=0{} & *=0{\mbox{\tiny\textbullet}} \ar[l] \ar@{-}[d] & \hspace{.8cm}\ar[l] & & & &
    *=0{\mbox{\tiny\textbullet}} \ar@{-}[d] \\
    & *=0{} & \hspace{0.1cm}\cdots\hspace{0.1cm} \ar[l] &&&\ar@{|->}[l]&
    *=0{\mbox{\tiny\textbullet}}
    \ar@{.}[d] \jumpdir{ E_1\lklammer{1.4cm}}{/:a(-45) +0.9cm/}
    \jumpdir{\rklammer{1.7cm}E}{/:a(78) +1.5cm/}\\
    && \hspace{0.1cm}\cdots\hspace{0.1cm} & *=0{} \ar[l]\ar@{-}[d]
    \jumpdir{G_{c-1}}{/:a(180) +0.2cm/} & & & *=0{\mbox{\tiny\textbullet}}
    \ar@{-}[d]\\
    & & \hspace{.8cm} & *=0{\mbox{\tiny\textbullet}} \ar@{-}[d] \ar[l] &
    *=0{} \ar[l] \ar@{-}[d]
    \jumpdir{N}{/:a(180) +0.2cm/} & & *=0{\mbox{\tiny\textbullet}} \ar@{-}[d]
    \jumpdir{\big\} E_c}{/:a(55) +0.5cm/} \\
    & & & *=0{} & *=0{} \ar[l] & & *=0{\mbox{\tiny\textbullet}}
  }
\end{equation}

For $c$ composable $1$-extensions $0 \leftarrow L_i \leftarrow G_i \leftarrow L_{i+1} \leftarrow 0$, with corresponding $1$-cocycles $\eta^{L_i}_{L_{i+1}}$, we define the set of admissible $c$-extension modules
\begin{eqnarray*}
 \lefteqn{\ExtMod(\eta^{L_0}_{L_1},\ldots,\eta^{L_{c-1}}_{L_c}):=} \\
 & & \{ E \mbox{ $c$-extension module} \mid \eta^{L_i}_{L_{i+1}} \mbox{ corresponds to } G_i(E) \mbox{ for all } i\}/\approx,
\end{eqnarray*}
where $E$ and $E'$ are equivalent ($E \approx E'$), if there exists an isomorphism $\alpha:E\to E'$, such that $\alpha$ induces isomorphisms from $E_i$ onto $E'_i$ for all $i$. For a $c$-extension module in $\ExtMod(\eta^{L_0}_{L_1},\ldots,\eta^{L_{c-1}}_{L_c})$ we occasionally use the terminology \emph{$(\eta^{L_0}_{L_1},\ldots,\eta^{L_{c-1}}_{L_c})$-extension module} to emphasize the dependency on the $1$-extension cocycles.

\begin{defn}[Rigid tuple of $1$-cocycles]
  We call a $c$-tuple $(\eta^{L_0}_{L_1},\ldots,\eta^{L_{c-1}}_{L_c})$ of {\sc Yoneda}-composable $1$-cocycles \emph{rigid}, if $|\ExtMod(\eta^{L_0}_{L_1},\ldots,\eta^{L_{c-1}}_{L_c})|=1$. In this case we also call the unique $(\eta^{L_0}_{L_1},\ldots,\eta^{L_{c-1}}_{L_c})$-extension module \emph{rigid}.
\end{defn}

In the Section~\ref{examples} we will provide examples for both rigid and non-rigid extension modules.

The natural problem that arises is to describe the set $\ExtMod(\eta^{L_0}_{L_1},\ldots,\eta^{L_{c-1}}_{L_c})$ and in particular to find necessary and sufficient conditions for it to be non-empty.

\section{The Case of $2$-Extension Modules}\label{case2}

Since for $c=1$ the notion of $1$-extensions and $1$-extension modules coincide, the first interesting case is $c=2$.

\[
\xymatrix@C=1.7cm@R=0.5cm{
	*=0{} \ar@{-}[d] \jumpdir{M}{/:a(180) +0.2cm/} & *=0{} \ar[l]
        \ar@{-}[d] \jumpdir{G_0}{/:a(180) +0.2cm/} &  &  \\
	*=0{} & *=0{\mbox{\tiny\textbullet}} \ar[l] \ar@{-}[d]
        \jumpdir{ L \cong \big\{ } {/:a(-55) +0.7cm/} & *=0{} \ar[l] \ar@{-}[d]
        \jumpdir{G_1}{/:a(180) +0.2cm/} \jumpdir{ \big\} L } {/:a(45) +0.5cm/}& \\
	 & *=0{} & *=0{\mbox{\tiny\textbullet}} \ar@{-}[d] \ar[l] &
         *=0{} \ar[l] \ar@{-}[d]
         \jumpdir{N}{/:a(180) +0.2cm/}\\
	 & & *=0{} & *=0{} \ar[l] \\
}
\]

Here we focus on the case of $2$-extensions. Let $G:0 \leftarrow M \leftarrow G_0 \leftarrow G_1 \leftarrow N \leftarrow 0$ with corresponding $2$-cocycle $\eta^G$ be the composite of $0 \leftarrow M \xleftarrow{\pi^{G_0}} G_0 \xleftarrow{\iota_{G_0}} L \leftarrow 0$ and $0 \leftarrow L \xleftarrow{\pi^{G_1}} G_1 \xleftarrow{\iota_{G_1}} N \leftarrow 0$ for $L:=L_1$. Let $\eta^M_L$ and $\eta^L_N$ denote the corresponding $1$-cocycles: $\eta^M_L\circ\eta^L_N=\eta^G$.

The short exact sequence $0 \leftarrow L \leftarrow G_1 \leftarrow N \leftarrow 0$ and the covariant functor $\Hom(M,-)$ give rise to the long exact sequence
\[
  \xymatrix@1{
    \cdots \Hom(M,L) \ar[r]^{\delta^0} & \Ext^1(M,N) \ar[r]^{{\iota_{G_1}}_*} & \Ext^1(M,G_1) \ar[r]^{{\pi^{G_1}}_*} &
    \Ext^1(M,L) \ar[r]^{\delta^1} & \Ext^2(M,N) \cdots
  }
\]
The existence of a $(\eta^M_L,\eta^L_N)$-extension module is equivalent to the existence of a $1$-cocycle $\eta^M_{G_1}\in\Ext^1(M,G_1)$ with ${\pi^{G_1}}_*(\eta^M_{G_1})=\eta^M_L$: \\
Each such element $\eta^M_{G_1}$ induces a $1$-extension  $\Yon(\eta^M_{G_1})=0 \leftarrow M \leftarrow E \leftarrow G_1 \leftarrow 0$ such that the induced $1$-extension $0 \leftarrow M\leftarrow E/N \leftarrow G_1/N \leftarrow 0$ corresponds to $\eta^M_L$. Thus $E$ is an extension module in $\ExtMod(\eta^M_L,\eta^L_N)$. \\
Due to the exactness of the above long exact sequence the existence of $\eta^M_{G_1}$ is in turn equivalent to $\delta^1(\eta^M_L)=0\in\Ext^2(M,N)$. But since $\delta^1(\eta^M_L) = \eta^M_L \circ \eta^L_N$ (cf.~Section~\ref{yoneda}) and $\eta^M_L \circ \eta^L_N=\eta^G$ one concludes:

\begin{axiom}[Existence Theorem]\label{the thm}
  $\ExtMod(\eta^M_L,\eta^L_N)\neq \emptyset$, if and only if the {\sc Yoneda} product $\eta^M_L\circ\eta^L_N\in\Ext^2(M,N)$ vanishes.
\end{axiom}

In other words, a $2$-extension $G:0 \leftarrow M \leftarrow G_0 \leftarrow G_1 \leftarrow N \leftarrow 0$ is induced by a $2$-extension module, if and only if $G=0 \in \YExt^2(M,N)$.

The above proof provides a surjection
\[
   \left({\pi^{G_1}}_*\right)^{-1}(\eta^M_L) \twoheadrightarrow \ExtMod(\eta^M_L,\eta^L_N)
\]
from the fiber of ${\pi^{G_1}}_*$ over $\eta^M_L$ (recall $G_1=\Yon(\eta^L_N)$) onto $\ExtMod(\eta^M_L,\eta^L_N)$.

Further, the long exact $\Ext$-sequence induces an action of $\Ext^1(M,N)$ on $\left({\pi^{G_1}}_*\right)^{-1}(\eta^M_L)$ given by
\[
  \begin{array}{rclcc}
   \Ext^1(M,N) &\times& \left({\pi^{G_1}}_*\right)^{-1}(\eta^M_L) & \to & \left({\pi^{G_1}}_*\right)^{-1}(\eta^M_L): \\
    \rule{0pt}{0.6cm}(\eta &,& \eta^M_{G_1}) & \mapsto & {\iota_{G_1}}_*(\eta) + \eta^M_{G_1} = \eta \circ \iota_{G_1} + \eta^M_{G_1}.
  \end{array}
\]
Because of the exactness of the above $\Ext$-sequence at $\Ext^1(M,G_1)$ this affine action is transitive. The kernel of the action is the kernel of ${\iota_{G_1}}_*$, which coincides due to the exactness at $\Ext^1(M,N)$ with the image of $\delta^0$:
\[
 \img \delta^0 = \Hom(M,L)\circ \eta^L_N := \{\phi^M_L \circ \eta^L_N \mid \phi^M_L\in\Hom(M,L) \},
\]
where the first equality was established in Section~\ref{yoneda}. Hence, the above action of $\Ext^1(M,N)$ turns the fiber $\left({\pi^{G_1}}_*\right)^{-1}(\eta^M_L)$ into a principal homogeneous space for
\[
  \Ext^1(M,N)/(\Hom(M,L)\circ \eta^L_N).
\]

\bigskip
Dually, the short exact sequence $0 \leftarrow M \leftarrow G_0 \leftarrow L \leftarrow 0$ and the contravariant functor $\Hom(-,N)$ give rise to the long exact sequence
\[
  \xymatrix@1{
    \cdots \Hom(L,N) \ar[r]^{\Delta^0} & \Ext^1(M,N) \ar[r]^{{\pi^{G_0}}^*} & \Ext^1(G_0,N) \ar[r]^{{\iota_{G_0}}^*} &
    \Ext^1(L,N) \ar[r]^{\Delta^1} & \Ext^2(M,N) \cdots
  }
\]
with $\Delta^1(\eta^L_N)=\eta^M_L\circ\eta^L_N$ (cf.~Section~\ref{yoneda}). This leads to a second surjection
\[
   \left({\iota_{G_0}}^*\right)^{-1}(\eta^L_N) \twoheadrightarrow \ExtMod(\eta^M_L,\eta^L_N).
\]
$\Ext^1(M,N)$ acts affinely on $\left({\iota_{G_0}}^*\right)^{-1}(\eta^L_N)$ via
\[
\begin{array}{rclcc}
   \Ext^1(M,N) &\times& \left({\iota_{G_0}}^*\right)^{-1}(\eta^L_N) & \to & \left({\iota_{G_0}}^*\right)^{-1}(\eta^L_N):   \\
   \rule{0pt}{0.6cm} (\eta &,& \eta^{G_0}_N) & \mapsto & {\pi^{G_0}}^*(\eta) + \eta^{G_0}_N = \pi^{G_0} \circ \eta + \eta^{G_0}_N,
  \end{array}
\]
with kernel of action being
\[
  \img \Delta^0 = \eta^M_L\circ\Hom(L,N).
\]
This turns the fiber $\left({\iota_{G_0}}^*\right)^{-1}(\eta^L_N)$ into a principal homogeneous space for
\[
  \Ext^1(M,N)/(\eta^M_L\circ\Hom(L,N)).
\]

The two surjections motivate the following theorem for which a ``coordinate'' proof is provided in Subsection~\ref{coord extmod2}.

\begin{axiom}[Modelling Theorem]\label{action}
  $\ExtMod(\eta^M_L,\eta^L_N)$ is a principal homogeneous space\footnote{Here we do not make any statement about the existence of a natural group structure for $\ExtMod(\eta^M_L,\eta^L_N)$, i.e.\  the existence of a naturally distinguished element in the ``affine'' space $\ExtMod(\eta^M_L,\eta^L_N)$ and its uniqueness. We leave this for future work.} for the abelian group
  \[
     \Ext^1(\eta^M_L,\eta^L_N) := \Ext^1(M,N)/(\Hom(M,L)\circ \eta^L_N+\eta^M_L\circ\Hom(L,N))
  \]
  and hence
  \[
    \begin{array}{lcrcl}
      \ExtMod(\eta^M_L,\eta^L_N)
      &\cong& \left({\pi^{G_1}}_*\right)^{-1}(\eta^M_L)  & / & (\eta^M_L\circ\Hom(L,N)) \\
      &\cong& \left({\iota_{G_0}}^*\right)^{-1}(\eta^L_N) & / & (\Hom(M,L)\circ \eta^L_N).
    \end{array}
  \]
  We call $\Ext^1(\eta^M_L,\eta^L_N)$ the first extension group of $1$-cocycles.
\end{axiom}

\begin{coro}
  A pair $(\eta^M_L,\eta^L_N)$ of $1$-cocycles is rigid, iff a $(\eta^M_L,\eta^L_N)$-extension modules exists and $\Ext^1(\eta^M_L,\eta^L_N)=0$.
\end{coro}

The argument leading to Theorem~\ref{the thm} expressed in ``coordinates'' will result in a linear inhomogenous system of equations (cf.~Subsection~\ref{coord extmod2}). In Subsections~\ref{compute yoneda} and \ref{triangular} we will see how Theorem~\ref{the thm} provides an alternative interpretation of this system that even reveals its triangular structure, which is extremely valuable for computations.

\section{The Higher $c$-Extension Modules ($c\geq 3$)}\label{3 and higher}

A necessary condition for the existence of higher $c$-extension modules follows immediately from Theorem~\ref{the thm}:
\begin{coro}\label{necessary}
  If $\ExtMod(\eta^{L_0}_{L_1},\ldots,\eta^{L_{c-1}}_{L_c})\neq \emptyset$ then the {\sc Yoneda} products $\eta^{L_{i-1}}_{L_i}\circ\eta^{L_i}_{L_{i+1}} \in \Ext^2(L_{i-1},L_{i+1})$ vanish for all $i=1,\ldots, c-1$.
\end{coro}

In Example~\ref{example not sufficient} we will provide an example showing that this condition is \emph{not} sufficient.

The following theorem is an obvious generalization of the argument preceeding Theorem~\ref{the thm}. It provides necessary and sufficient conditions for the existence of a $c$-extension module with given $(\eta^{L_0}_{L_1},\ldots,\eta^{L_{c-1}}_{L_c})$ by using Theorem~\ref{the thm} as the induction step.

\begin{coro}\label{coro}
  The set $\ExtMod(\eta^{L_0}_{L_1},\ldots,\eta^{L_{c-1}}_{L_c})$ is non-empty, if and only if
  \begin{enumerate}
  \item $\ExtMod(\eta^{L_1}_{L_2},\ldots,\eta^{L_{c-1}}_{L_c})\neq \emptyset$ \emph{and}
  \item there exists a $(c-1)$-extension module $E_1 \in \ExtMod(\eta^{L_1}_{L_2},\ldots,\eta^{L_{c-1}}_{L_c})$, such that the {\sc Yoneda} product \[
    \eta^{L_0}_{L_1}\circ \eta^{L_1}_{E_2} = 0 \in \Ext^2(L_0,E_2),
  \]
  where the $1$-cocycle $\eta^{L_1}_{E_2}$ is induced by $E_1$.
  \end{enumerate}
\end{coro}

At the end of Subsection~\ref{coord extmod3} we will use Theorem~\ref{the thm} to provide a simple example showing that there may very well exist a $(c-1)$-extension module $E_1 \in \ExtMod(\eta^{L_1}_{L_2},\ldots,\eta^{L_{c-1}}_{L_c})$ with $\eta^{L_0}_{L_1}\circ \eta^{L_1}_{E_2} \neq 0 \in \Ext^2(L_0,E_2)$, i.e.\  that does not lead to a $c$-extension module, whereas a different choice of $E_1$ does. This narrows the range of applicability of Corollary~\ref{coro} considerably for $c>3$. This limitation will be explained in Subsection~\ref{coord 4 and higher}. The special case $c=3$ will be discussed in Subsection~\ref{coord extmod3}.

\section{The ``Coordinate'' Description}\label{coord}

By a ``coordinate description'' of a module $M$ we simply mean a \emph{finite free presentation} given by a matrix $d_1 \in D^{p\times q}$, where $d_1$ is viewed as a morphism of free modules $d_1:D^{1\times p} \to D^{1\times q}$. $M$ is then the cokernel of $d_1$. $d_1$ is therefore called \emph{presentation matrix} or \emph{matrix of relations}. This yields the beginning of a free resolution
\[
   0 \leftarrow M \leftarrow F_0 \xleftarrow{d_1} F_1,
\]
where $F_0:=D^{1\times q}$ and $F_1:=D^{1\times p}$.

In most of the arguments used below the modules need not be finitely generated and one can replace the word matrix by morphism or matrix of morphisms. In particular, the proof of Theorem~\ref{action} given in Subsection~\ref{coord extmod2} applies without the restriction of being finitely presented.

\subsection{$\Ext^1$ in ``coordinates''}\label{coord ext1} In this subsection we recall well-known facts about $\Ext^1$. The {\sc Yoneda} correspondence between $c$-cocycles and $c$-extension is for $c=1$ summed up in the diagram (cf.~Appendix~\ref{equivalence})
\[
  \xymatrix{
    0 & \ar[l] M \ar@{=}[d] & \ar[l] F_0 \ar[d]^{\eta_0} & \ar[l]_{d_1}
  K_1 \ar[d]^{\eta} & \ar[l]   0 \\
    0 & \ar[l] M & \ar[l]_{\pi} E & \ar[l]_{\iota} N & \ar[l] 0,
  }
\]
which shows how to compute $\eta$ by lifting the identity $M\xrightarrow{\id}M$ twice. Conversely, $M$ is in the above diagram the \emph{pushout} of $N\xleftarrow{\eta} K_1
\xrightarrow{d_1} F_0$, i.e.\  the cokernel of
\[
  K_1 \xrightarrow{\left(\begin{array}{cc} d_1 & \eta
\end{array}\right)} F_0 \oplus N \xrightarrow{\left(\begin{array}{c}
\eta_0 \\ -\iota \end{array}\right)} E \to 0.
\]

Following the notational convention in \cite{BR}: If we write $M$ for the cokernel of the relation matrix $d_1$ then we write $\tM$ for $d_1$. The module $E$ in the above sequence is then the cokernel of the matrix (recall, we use the row convention)
\[
  \tE := \left(\begin{array}{c|c}
                  \tM & \eta \\ \hline
                  0 & \tN 
                \end{array} \right).
\]
The upper row is the morphism $\left(\begin{array}{cc} d_1 & \eta
\end{array}\right)$, whereas the second row is a presentation matrix for the module $F_0\oplus N$.

Now we will make an attempt to derive the defining properties of a $1$-cocycle (cf.~Appendix~\ref{satellites})
\begin{equation}\tag{$\Ext^1$}\label{ext1}
   \eta \in \Ext^1(M,N):=\frac{\{\eta: F_1 \to N \mid  0=d_2\eta\}}
   {\{ d_1\phi \mid \phi: F_0 \to N\}}
\end{equation}
in an elementary way: \\
Recall, an extension of $N$ by $M$ is described by a module $E$ of which $N$ is a submodule and $M=E/N$ is the factor module modulo $N$. Now we want to explicitly construct such an $E$ as the cokernel of a matrix $\tE$. We want the standard basis row vectors of the form $\left(\begin{array}{ccc|ccccc} 0 & \cdots & 0 & 0 & \cdots& 1 & \cdots & 0\end{array}\right)$ to be the representatives of generators of $N$ in the cokernel of $\tE$, so the lower part of the matrix can now be set to $\left(\begin{array}{c|c} 0 & \tN \end{array}\right)$, where $\tN$ is a relation matrix for $N$. Now computing modulo these vectors $\left(\begin{array}{c|c} 0 & 1 \end{array}\right)$ states that $E/N$ is presented by the left hand side of the matrix $\tE$. Putting a presentation matrix $\tM$ with cokernel $M$ in the upper left corner (above the zeros) thus leads to a factor module isomorphic to $M$. Now $\tE$ has the form $\left(\begin{smallmatrix} \tM & \eta \\ 0 & \tN            \end{smallmatrix} \right)$ and the conditions on $\eta$ remain to be determined: $N$ is described by \emph{all} the relations among the vectors $\left(\begin{array}{c|c} 0 & 1 \end{array}\right)$ and to isolate them one needs the {\em most general} row operation matrix, which applied to the upper rows leads to a zero matrix on the left hand side. This is precisely the \emph{first syzygies matrix} $d_2:F_2\to F_1$ satisfying $d_2 \tM = d_2 d_1 = 0$:
\[
   \left(\begin{array}{c|c} d_2 & 0 \end{array}\right)\tE = \left(\begin{array}{c|c}
                  0 & d_2\eta \\ \hline
                  0 & \tN 
                \end{array} \right).
\]
Thus $d_2\eta$ must not introduce new relations to $\tN$, which means $d_2\eta = 0$ modulo the relations $\tN$. In other words $d_2\eta=0$ as a morphism $F_2\to N$. This gives back the numerator of (\ref{ext1}). To explain the denominator we proceed as follows: We now consider the most general \emph{invertible} row and column operations on $\tE$ preserving the above situation, i.e.\  preserving the submodule $N$ with factor module $M$. These operations lead to a congruent extension (cf.~(\ref{cong})~in~Appendix~\ref{composition}). In particular we want to preserve the \emph{upper triangular structure} of the matrix $\tE$, i.e.\  zeros in the lower left corner. Without loss of generality (neglecting possible ``coordinate changes'' of $M$ and $N$) one can even consider only those operations that leave the submatrices $\tM$ and $\tN$ in $\tE$ fixed. This leaves us out with the following two possibilities:
\[
  \xymatrix{
    \mbox{$\left(\begin{array}{c|c}
       \tM & \eta \\ \hline
       0 & \tN 
    \end{array}
  \right)$}\jumpdir{\begin{turn}{-90}$\curvearrowleft$\end{turn}} {/:a(-88)
  +1cm/}\jumpdir{\text{\footnotesize \mbox{$\ +\ \chi\ \cdot$}}}{/:a(-90) +1.6cm/}}
  \quad\mbox{ and }\quad
  \xymatrix{
    \mbox{$\left(\begin{array}{c|c}
       \tM & \eta \\ \hline
       0 & \tN 
    \end{array} \right).$}\jumpdir{\stackrel{+\  \cdot\  \phi}{\curvearrowright}}{/:a(0) +0.8cm/}
  }
\]
The row operation only replaces $\eta$ by $\eta + \chi \tN$, which does not change $\eta$ considered as a morphism $F_1\to N$. The column operation replaces $\eta$ by $\eta + \tM \phi=\eta + d_1\phi$, giving back the denominator in (\ref{ext1}).

\bigskip
Since we are working in coordinates we actually have to destinguish between cocycles and the matrices representing them. In the next subsection a strict distinction in the notation will be unavoidable. So from now on we denote a matrix representative of a cocycle $\eta$ by $\tilde\eta$. Conversely, the cocycle represented by a matrix $\tilde\eta$ is denoted by $\eta=[\tilde{\eta}]$.

\subsection{$\ExtMod(\eta^M_L,\eta^L_N)$ in ``coordinates'' ($c=2$)}\label{coord extmod2} In this subsection we provide a coordinate description and proof of Theorem~\ref{action}.

Let $M$, $L$ and $N$ be as in  Section \ref{case2}. Following the line of argument of the previous subsection one can assume that a matrix of relations of a $2$-extension module in $\ExtMod(\eta^M_L,\eta^L_N)$ has the upper triangular form
\[
   \tE := \left(\begin{array}{ccc}
                  \tM & \tilde\eta^M_L & \eta  \\
                  \cdot & \tL & \tilde\eta^L_N \\
                  \cdot & \cdot & \tN 
          \end{array}\right).
\]
The ``coordinate'' description of $\ExtMod(\eta^M_L,\eta^L_N)$ boils down to
classifying the admissible $\eta$'s. For the cokernel of $\tE$ to be in
$\ExtMod(\eta^M_L,\eta^L_N)$ it is necessary and sufficient for
$\left(\begin{array}{cc} \tilde\eta^M_L & \eta \end{array}\right)$ to describe a
cocycle in
\[
\Ext^1(M,G_1)=\Ext^1(\coker(\tM),\coker\left(\begin{array}{cc}\tL & \tilde\eta^L_N
    \\ 0 & \tN \end{array}\right))
\]
i.e.
\begin{equation*}
  d_2^M \left(\begin{array}{cc} \tilde\eta^M_L & \eta \end{array}\right) =
  \left(\begin{array}{c|c}X_1 & X_2\end{array}\right)
  \left(\begin{array}{cc}\tL & \tilde\eta^L_N \\ 0 & \tN \end{array}\right),
\end{equation*}
where $d_2^M$ is the first syzygies matrix of $M=\coker(\tM)$. This leads to the \emph{\textbf{linear}
inhomogenous} system of equations with coefficients in $D$:
\begin{equation}\tag{$\ExtMod^2$}\label{extmod2}
\begin{array}{rlcl}
  d_2^M & \tilde\eta^M_L & = & \enbox{0}{-0.15}{0.7}{0.6}~X_1 \  \tL \\
  d_2^M & \enbox{-0.05}{-0.15}{0.6}{0.6}~\eta \rule{0pt}{0.6cm} & = & \enbox{0}{-0.15}{0.7}{0.6}~X_1\  \tilde\eta^L_N + \enbox{0}{-0.15}{0.7}{0.6}~X_2 \  \tN
\end{array}
\end{equation}
Equivalently, $\left(\begin{array}{cc} \eta \\ \tilde\eta^L_N \end{array}\right)$
must describe a cocycle in
\[
\Ext^1(G_0,N)=\Ext^1(\coker\left(\begin{array}{cc}\tM & \tilde\eta^M_L \\ 0 & \tL
  \end{array}\right),\coker \tN),
\]
which leads to the \emph{same} system of equations.

\bigskip
Again, without loss of generality one can assume $\tM,\tL$ and $\tN$ fixed. First we define the set of all matrix representatives of the fixed pair of $1$-cocycles $(\eta^M_L,\eta^L_N)$
\[
  B(\eta^M_L,\eta^L_N) := \left\{(\tilde{\eta}^M_L,\tilde{\eta}^L_N)\mid [\tilde{\eta}^M_L]=\eta^M_L,\
   [\tilde{\eta}^L_N]=\eta^L_N\right\}.
\]
Second we consider the disjoint union
\begin{eqnarray*}
  \lefteqn{\widetilde{\ExtMod}(\eta^M_L,\eta^L_N) := }\\
   & & \rule{0pt}{0.8cm}\coprod_{ (\tilde{\eta}^M_L,\tilde{\eta}^L_N)\in B(\eta^M_L,\eta^L_N) }  \{ \eta \mid 
   \mbox{ the matrix } \eta \mbox{ satisfies (\ref{extmod2}) for the matrices }
   \tilde{\eta}^M_L, \tilde{\eta}^L_N \}
\end{eqnarray*}
together with the obvious projection $\pi:\widetilde{\ExtMod}(\eta^M_L,\eta^L_N) \twoheadrightarrow B(\eta^M_L,\eta^L_N)$. We want to identify elements of $\widetilde{\ExtMod}(\eta^M_L,\eta^L_N)$ with relation matrices of the form
\[
   \tE = \left(\begin{array}{ccc}
                  \tM & \tilde{\eta}^M_L & \eta  \\
                  \cdot & \tL & \tilde{\eta}^L_N \\
                  \cdot & \cdot & \tN 
          \end{array}\right).
\]
This also emphasizes the dependency of the two sets just introduced on the choice of the presentation matrices $\tM$, $\tL$, and $\tN$. Further we introduce the double-unipotent group
\[
 U := \left\{
        \left(
          \begin{array}{ccc}
            1 & \kappa & \chi \\
            0 & 1 & \nu \\
            0 & 0 & 1
          \end{array}
        \right)
      \right\}
      \times
      \left\{
        \left(
          \begin{array}{ccc}
            1 & \mu & \phi \\
            0 & 1 & \lambda \\
            0 & 0 & 1
          \end{array}
        \right)
      \right\},
\]
together with its action on $\widetilde{\ExtMod}(\eta^M_L,\eta^L_N)$ given by
\[
  \begin{array}{cclcc}
    U & \times & \widetilde{\ExtMod}(\eta^M_L,\eta^L_N) & \to & \widetilde{\ExtMod}(\eta^M_L,\eta^L_N): \\
    ((r,c) & , & \tE) & \mapsto & r \, \tE \, c^{-1}.
  \end{array}
\]
Since $U$ is the biggest group fixing the diagonal and preserving the triangular structure of $\tE$, there is a 1-1-correspondence between the set $\ExtMod(\eta^M_L,\eta^L_N)$ and the global quotient\footnote{This argument generalizes to $c>2$ in the obvious way.} $\widetilde{\ExtMod}(\eta^M_L,\eta^L_N)/U$.

Before we proceed we illustrate the action by the following row and column operations:
\[
  \xymatrix{
    \mbox{$\left(\begin{array}{ccc}
                  \tM & \tilde\eta^M_L & \eta  \\
                  \cdot & \tL & \tilde\eta^L_N \\
                  \cdot & \cdot & \tN 
          \end{array}\right)$}\jumpdir{\begin{tikzpicture}
  \draw[<-] (0,0.9) .. controls (0.5,0.72) and (0.4,0.1) .. (0,0);
\end{tikzpicture}}{/:a(-90) +1.7cm/}
      \jumpdir{\text{\footnotesize \mbox{$\  \  +\  \chi\  \cdot$}}}{/:a(-90)
        +2.2cm/}
      \jumpdir{\stackrel{-\  \cdot\
          \phi}  {\begin{tikzpicture}
  \draw[->] (0,0) .. controls (0.2,0.5) and (1.2,0.5) .. (1.4,0);
\end{tikzpicture}}}{/:a(5) +1.2cm/}
  }
  \xymatrix{
    \mbox{$,\left(\begin{array}{ccc}
                  \tM & \tilde\eta^M_L & \eta  \\
                  \cdot & \tL & \tilde\eta^L_N \\
                  \cdot & \cdot & \tN 
          \end{array}\right)$}
     \jumpdir{\begin{tikzpicture}
  \draw[<-] (0,0.5) .. controls (0.4,0.4) and (0.4,0.1) .. (0,0);
\end{tikzpicture}}{/:a(-83) +1.8cm/}
      \jumpdir{\text{\footnotesize \mbox{$\  \  +\  \kappa\  \cdot$}}}{/:a(-87)
        +2.3cm/}\jumpdir{\stackrel{-\  \cdot\
          \lambda}{\begin{tikzpicture}
  \draw[->] (0,0) .. controls (0.1,0.4) and (0.7,0.4) .. (0.8,0);
\end{tikzpicture}}}{/:a(-17) +1.3cm/}
  }
  \xymatrix{
    \mbox{$,\left(\begin{array}{ccc}
                  \tM & \tilde\eta^M_L & \eta  \\
                  \cdot & \tL & \tilde\eta^L_N \\
                  \cdot & \cdot & \tN 
          \end{array}\right)\hspace{1.3cm}.$}\jumpdir{\begin{tikzpicture}
  \draw[<-] (0,0.5) .. controls (0.4,0.4) and (0.4,0.1) .. (0,0);
\end{tikzpicture}}{/:a(-104) +1.1cm/}
      \jumpdir{\text{\footnotesize \mbox{$\  \  +\  \nu\  \cdot$}}}{/:a(-100)
        +1.6cm/}
      \jumpdir{\stackrel{-\  \cdot\
          \mu}{\begin{tikzpicture}
  \draw[->] (0,0) .. controls (0.1,0.4) and (0.7,0.4) .. (0.8,0);
\end{tikzpicture}}}{/:a(42) +1.6cm/}
  }
\]
Since $U$ acts on the set of fibers of $\pi$ it acts equivalently on the set $B(\eta^M_L,\eta^L_N)$ and this action turns out to be transitive: $\tilde\eta^M_L \to \tilde\eta^M_L+\kappa\,\tL-\tM\,\mu$ and $\tilde\eta^L_N \to \tilde\eta^L_N + \nu\,\tN-\tL\,\lambda$.

The orbit space $\widetilde{\ExtMod}(\eta^M_L,\eta^L_N)/U$ is due to the transitivity of the induced action on $B(\eta^M_L,\eta^L_N)$ naturally bijective to the global quotient
\[
\pi^{-1}((\tilde{\eta}^M_L,\tilde{\eta}^L_N))/\mathrm{Stab}_{(\tilde{\eta}^M_L,\tilde{\eta}^L_N)}(U)
\]
for an arbitrary but fixed pair of matrices $(\tilde{\eta}^M_L,\tilde{\eta}^L_N)\in B(\eta^M_L,\eta^L_N)$.

$\mathrm{Stab}_{(\tilde{\eta}^M_L,\tilde{\eta}^L_N)}(U)$ is thus the largest subgroup which fixes the secondary diagonal and we conclude that
\[
  \mathrm{Stab}_{(\tilde{\eta}^M_L,\tilde{\eta}^L_N)}(U) =
  \left\{(
    \left(\begin{array}{ccc}
            1 & \kappa & \chi \\
            0 & 1 & \nu \\
            0 & 0 & 1
          \end{array}
        \right),
        \left(
          \begin{array}{ccc}
            1 & \mu & \phi \\
            0 & 1 & \lambda \\
            0 & 0 & 1
          \end{array}\right) )
    \mid \kappa\,\tL-\tM\,\mu = 0 \mbox{ and } \nu\,\tN-\tL\,\lambda = 0 \right\}.
\]
There are no conditions on the matrices $\phi$ and $\chi$, while the two specified conditions are interpreted as follows: \\
$\tL\,\lambda=\nu\,\tN$ states that $\lambda$ defines a morphism in $\Hom(L,N)$ and hence\footnote{For the computation of the \textsc{Yoneda} product cf.~Appendix~\ref{lift}.} 
\[
\tilde\eta^M_L\,\lambda \in \eta^M_L\circ \Hom(L,N).
\]
$\kappa\,\tL=\tM\,\mu$ states that $\mu$ is a morphism in $\Hom(M,L)$ and $\kappa$ is a lift of $\mu$ in the diagram
\[
\xymatrix@C=1.2cm{
  0 & M \ar^{\mu}[d] \ar[l] & F_0 \ar_{I_{\mathrm{rank}(F_0)}}[l]  \ar@{.>}^{\mu}[d] & F_1 \ar_{d_1^M = \tM}[l]
  \ar@{.>}_{\kappa}[d] \ar@{-->}@/^1.5pc/^{\text{$\mu\circ \eta^L_N$}}[dd]\\
  0 & L \ar[l] & F'_0 \ar^{I_{\mathrm{rank}(F_0')}}[l] & F'_1 \ar^{d_1^L =
    \tL}[l] \ar_{\tilde\eta^L_N}[d]\\
  & & & \tN
}
\]
thus 
\[
\kappa\,\tilde\eta^L_N \in \Hom(M,L) \circ \eta^L_N.
\]

Hence, the stabilizer leads to the four operations
\[
  \begin{array}{rclcl}
  \eta & \to & \eta+\chi\,\tN & & (\chi \mbox{ arbitrary}) \\
  \eta & \to & \eta-\tM\,\phi & \quad & (\phi \mbox{ arbitrary}) \\
  \eta & \to & \eta+\kappa\ \tilde\eta^L_N\rule{0cm}{0.5cm} & & (\kappa \mbox{ lift}) \\
  \eta & \to & \eta-\tilde\eta^M_L\,\lambda\rule{0pt}{0.5cm} & & (\lambda \mbox{ morphism})
  \end{array}
\]
defining the \emph{congruence relation} between different $\eta$'s that represent the same element in $\ExtMod(\eta^M_L,\eta^L_N)$:
\[
  \eta \approx \eta' \quad \mbox{ iff } \quad \eta-\eta' = \chi\,\tN-\tM\,\phi+\kappa\,\tilde\eta^L_N-\tilde\eta^M_L\,\lambda, \quad \mbox{ where } \phi,\chi,\kappa,\lambda \mbox{ as above}.
\]
Summing up, one can identify $\ExtMod(\eta^M_L,\eta^L_N)$ with the set of all $\eta$'s satisfying the above system of equations (\ref{extmod2}) \emph{modulo} the congruence relation $\approx$.

\begin{proof}[Proof of Theorem~\ref{action}]
It is now easy to see that
\[
  \eta \mapsto \eta+\tilde\eta^M_N
\]
defines a natural action of $\Ext^1(M,N)$ on the set of all $\eta$ satisfying (\ref{extmod2}) with $\Hom(M,L)\circ \eta^L_N+\eta^M_L\circ\Hom(L,N)$ being the largest subgroup that acts trivially:\\
A cocycle in $\Ext^1(M,N)$ is represented by a morphism $\eta^M_N\in\Hom(K_1,N)$ which is in turn represented by a matrix $\tilde\eta^M_N$ that fullfills the equation $d_2^M \tilde\eta^M_N = X_3 \tN$ for some matrix $X_3$ (compare with the numerator of (\ref{ext1})). Thus $\eta+\tilde\eta^M_N$ still satisfies the above equations (\ref{extmod2}), with $X_2$ replaced by $X_2+X_3$. Moreover, a second matrix $\bar\eta^M_N$ representing the same morphism $\eta^M_N\in\Hom(K_1,N)$ differs from $\tilde\eta^M_N$ by a matrix of the form $\chi\,\tN$ and hence $\eta+\bar\eta^M_N \approx \eta+\tilde\eta^M_N$.  \\
For $\Ext^1(M,N)$ to act we need to verify that the subgroup $d_1\Hom(F_0,N)$ acts trivially (compare with the denominator of (\ref{ext1})). But since $d_1 = \tM$, the congruence relation asserts that $\eta \approx \eta+\tM\phi$ and the action of this subgroup is indeed trivial. \\
The remaining two operations $\eta \to \eta+\kappa\ \tilde\eta^L_N$ and $\eta \to \eta-\tilde\eta^M_L\,\lambda$ of the congruence relation state that the subgroup $\Hom(M,L)\circ \eta^L_N+\eta^M_L\circ\Hom(L,N)$ coincides with the kernel of the action of $\Ext^1(M,N)$ on $\ExtMod(\eta^M_L,\eta^L_N)$.

To complete the proof of Theorem~\ref{action} we still need to see that $\Ext^1(M,N)$ acts transitively on $\ExtMod(\eta^M_L,\eta^L_N)$. To this end let $\bar{\eta}, \bar{X}_1,\bar{X}_2$ be a second solution of the system (\ref{extmod2}). Then $(X_1 - \bar{X}_1)\tL= (d_2^M\tilde\eta^M_L-d_2^M\tilde\eta^M_L)=0$ and hence there exists a $Y$ such that $X_1-\bar{X}_1 = Y d_2^L$, where $d_2^L$ is the first syzygies matrix of $L=\coker(\tL)$. It follows that $d_2^M(\eta-\bar{\eta})=(X_1-\bar{X}_1)\tilde\eta^L_N+(X_2-\bar{X}_2)\tN = Y d_2^L \tilde\eta^L_N + (X_2-\bar{X}_2)\tN$ and since $d_2^L\tilde\eta^L_N=0 \mod \tN$ we conclude $d_2^M(\eta-\bar{\eta})=0 \mod \tN$. This means $\eta-\bar{\eta}$ is indeed a matrix representing a cocycle.
\end{proof}

\bigskip
Refining the congruence relations $\approx$ in two different ways and fixing a pair $(\tilde\eta^M_L,\tilde\eta^L_N)$ with $[\tilde\eta^M_L]\in\Ext^1(M,L)$ and $[\tilde\eta^L_N]\in\Ext^1(L,N)$ we recover the two fibers from Section~\ref{case2}: The set of all $\left(\begin{array}{cc} \tilde\eta^M_L & \eta \end{array}\right)$ interpreted as $1$-cocycles in $\Ext^1(M,G_1)$ with $\eta$ satisfying (\ref{extmod2}) gives back $\left({\pi^{G_1}}_*\right)^{-1}(\eta^M_L)$. This means that we refine the congruence relation by dropping the operation $\eta \to \eta-\tilde\eta^M_L\,\lambda$ (or equivalently restricting the action to the subgroup of $U$ defined by $\lambda=0$ and $\nu=0$), which would in general alter $[\left(\begin{array}{cc} \tilde\eta^M_L & \eta \end{array}\right)]\in\Ext^1(M,G_1)$. Similarly with the set of all $[\left(\begin{array}{cc} \eta \\ \tilde\eta^L_N \end{array}\right)]\in\Ext^1(G_0,N)$ and $\left({\iota_{G_0}}^*\right)^{-1}(\eta^L_N)$, dropping the operation $\eta \to \eta+\kappa\ \tilde\eta^L_N$. We illustrate the two situations in the matrix $\tE$:
\[
  \left(\begin{array}{ccc}
                  \tM & \multicolumn{2}{c}{\boxed{\begin{array}{cc} \tilde\eta^M_L\hspace{5pt} &
                        \eta\hspace{7pt} \end{array}}}  \\ 
                  \cdot & \  \  \  \tL & \tilde\eta^L_N\rule{0cm}{0.5cm} \\
                  \cdot & \  \  \  \cdot & \tN\rule{0cm}{0.5cm} 
          \end{array}\right)
\  \  \mbox{ and } \  \
\left(\begin{array}{ccc}
                  \tM\rule{0cm}{0.5cm} & \tilde\eta^M_L & \multirow{2}{*}{\boxed{\begin{array}{c}
                        \eta \\ \tilde\eta^L_N\rule{0cm}{0.5cm} \end{array}}} \\
                  \cdot & \tL\rule{0cm}{0.6cm}  \\
                  \cdot & \cdot & \tN\rule{0cm}{0.5cm} 
          \end{array}\right).
\]

\subsection{How to compute the {\sc Yoneda} product of two $1$-cocycles}\label{compute yoneda}
As explained in Appendix~\ref{lift}, the {\sc Yoneda} product of two $1$-cocycles $\eta^M_L \in \Ext^1(M,L)$ and $\eta^L_N \in \Ext^1(L,N)$ can be computed by lifting $\eta^M_L$ to $X_1$ in the diagram:
\[
\xymatrix@C=1.5cm{
  & & F_1 \ar[dl]_{\eta^M_L} \ar@{.>}^{\eta^M_L}[d] & F_2 \ar_{d_2^M}[l]
  \ar@{.>}_{\encircle{0.32}{0.08}{0.5}~X_1}[d] \ar@{-->}@/^1.5pc/^{\text{\fbox{$\eta^M_L\circ \eta^L_N$}}}[dd]\\
  0 & L \ar[l] & F'_0 \ar^{I_{\mathrm{rank}(F_0')}}[l] & F'_1 \ar^{d_1^L =
    \tL}[l] \ar_{\eta^L_N}[d]\\
  & & & \tN
}
\]
This means we have to solve the $D$-linear inhomogenous equation
\begin{equation*}
\begin{array}{rlcl}
  d_2^M & \eta^M_L & = & \encircle{0.375}{0.1}{0.7}~X_1\  \tL
\end{array}
\end{equation*}
and then simply compute the product $X_1\eta^L_N$. This is a matrix representing the {\sc Yoneda} product of the two cocycles:
\[
  \boxed{\  \rule[-0.2cm]{0pt}{0.7cm}\eta^M_L\circ \eta^L_N = X_1\eta^L_N.\  }
\]

\subsection{The system (\ref{extmod2}) is triangular}\label{triangular} 
Now we want to study the solvability of the system (\ref{extmod2}):
\[
\begin{array}{rlcl}
  d_2^M & \eta^M_L & = & \encircle{0.375}{0.1}{0.7}~X_1\  \tL \\
  d_2^M & \enbox{-0.05}{-0.15}{0.6}{0.6}~\eta & = & \enboxdash{0}{-0.15}{0.7}{0.6}~X_1\  \eta^L_N + \enbox{0}{-0.15}{0.7}{0.6}~X_2 \  \tN \rule{0pt}{0.7cm}
\end{array}
\]
by trying to describe its compatibility conditions. Using the insight of Subsection~\ref{compute yoneda}, we can reinterpret the system in the following way: The upper equation defines $X_1$ which always exists as the matrix of the lift (cf.~Subsection~\ref{compute yoneda}), \emph{even though} this equation is inhomogenous. Furthermore, for \emph{any} solution $X_1$, the product matrix $X_1\eta^L_N$ is a representative of the $2$-cocycle $\eta^M_L\circ \eta^L_N$ (cf.~Subsection~\ref{compute yoneda}). The lower equation precisely states that the {\sc Yoneda} product $\eta^M_L\circ \eta^L_N$ is zero as a $2$-cocycle, and hence the solvability of the system only depends on the {\sc Yoneda} product and \emph{not} on the choice of $X_1$.

\bigskip

The above discussion ``coordinatizes'' the proof of Theorem~\ref{the thm}, since it shows that the vanishing of the {\sc Yoneda} product is expressed by the solvability of the system. Hence, one has to admit, that at first glance nothing is won if one considers the vanishing of the {\sc Yoneda} product as the compatibility condition of the system\footnote{As if one would say, the system is solvable, if and only if the system is solvable.}. But in the course of showing this we discovered the system to be triangular\footnote{Recall that over a ring a triangular shape does not necessarily imply triangular structure, cf.~Appendix~\ref{appendix triangular}.} with respect to the shape given above, i.e.\  the system is successively solvable by first solving the upper equation and then the lower. This is of considerable computational value. Moreover, this triangular structure will play a decisive role for the case $c=3$ in Subsection~\ref{coord extmod3}.  \qed

\bigskip

It is important to note that the our notion of ``triangular system'' over a ring is well defined only with respect to a given triangular shape and we want to emphasize that our statement only applies to the triangular shape given above.

\subsection{$\ExtMod(\eta^M_L,\eta^L_K,\eta^K_N)$ in ``coordinates'' ($c=3$)}\label{coord extmod3}

In Corollary~\ref{coro} we gave an inductive condition for the existence of a $c$-extension module. This translates for $c=3$ into ``coordinates'' as follows: In the relation matrix
\[
\left(\begin{array}{cccc}
    \ \tM\  & \boxed{\eta^M_L}\rule[-0.2cm]{0cm}{0.6cm} & \eta^M_{K} & \eta^M_N \\
     \cdot  &  \tL &   \multicolumn{2}{c} {\boxed{\begin{array}{ccc} \eta^{L}_{K} &  \eta^{L}_N \rule[-0.2cm]{0cm}{0.7cm} \end{array}}}\\
     \cdot \rule[-0.1cm]{0cm}{0.6cm} &  \cdot &  \tK &  \eta^{K}_N \\
     \cdot &  \cdot &  \cdot &  \tN\rule{0cm}{0.5cm}
\end{array}\right)
\]
the submatrix $\left(\begin{array}{cc} \eta^{L}_{K} & \eta^{L}_N\end{array}\right)$ must be a $1$-cocycle in
\[
  \Ext^1(L,E_2)=\Ext^1(\coker(\tL),\coker\left(\begin{array}{cc}\tK & \eta^K_N \\ 0 & \tN \end{array}\right)) \quad  (\ast)
\]
and its {\sc Yoneda} product $\eta^M_L\circ\left(\begin{array}{cc} \eta^{L}_{K} & \eta^{L}_N\end{array}\right)$ with $\eta^M_L$ must vanish ($\ast\ast$). This is summed up in the quadratic system of equations (cf.~Subsections~\ref{coord ext1} and \ref{compute yoneda}):
\begin{flushleft}
$\begin{array}{lrcl}
    (\ast) \hspace{2cm}& d_2^{L} \left(\begin{array}{cc} \eta^{L}_{K} &
    \eta^{L}_N\end{array}\right) & = & \left(\begin{array}{c|c}Y_1 &
    Y_2\end{array}\right) \left(\begin{array}{cc}\tK & \eta^{K}_N \\ 0 &
     \tN \end{array}\right)\rule[-0.5cm]{0cm}{0.7cm}\\
  \cline{2-4} \setlength{\unitlength}{1cm}\begin{picture}(0,0)\put(2.71,0.625){\rule{11.77cm}{0.5pt}}\end{picture} 
  \setlength{\unitlength}{1cm}\begin{picture}(0,0)\put(0,-0.3){$(\ast\ast)$}\end{picture}& d_2^M \eta^M_L & = & X_1 \tL\rule{0cm}{0.7cm}\\
& d_2^{M} \left(\begin{array}{cc} \eta^M_{K} &
    \eta^M_N\end{array}\right) & = & X_1 \left(\begin{array}{cc}\eta^{L}_{K} & \eta^{L}_N \end{array}\right) + \left(\begin{array}{c|c}
    X_2 & X_3 \end{array}\right) \left(\begin{array}{cc}
      \tK & \eta^{K}_N \\ 0 & \tN \end{array}\right).
\end{array}$
\end{flushleft}
Again, as in Subsection~\ref{triangular}, the middle equation is always solvable and independent from the rest. Therefore, the rest of the quadratic system is in fact \emph{inhomogenous \textbf{linear}}:
\[
\begin{array}{r@{\hspace{0.2cm}}lcr@{\hspace{0.2cm}}lcr@{\hspace{0.2cm}}lcr@{\hspace{0.2cm}}l}
  d_2^{L}\rule{0cm}{0.7cm} & \eta^{L}_{K}    & = &
  \enbox{-0.2}{-0.2}{0.7}{0.7} Y_1 & \tK\\
  d_2^{L}\rule[-0.3cm]{0cm}{1cm} & \enbox{-0.1}{-0.2}{0.7}{0.7} \eta^{L}_N & = &
  \enbox{-0.2}{-0.2}{0.7}{0.7} Y_1 & \eta^{K}_N & + &
  \enbox{-0.2}{-0.2}{0.7}{0.7} Y_2 & \tN \\
  \hline \hline
  d_2^M\rule[-0.3cm]{0cm}{0.9cm}     & \eta^M_{L}        & = &
  \encircle{0.375}{0.1}{0.7}~X_1 & \tL\\
  \hline
  d_2^M\rule{0cm}{0.7cm}     & \enbox{-0.1}{-0.2}{0.7}{0.7} \eta^M_{K} & = &
  \enboxdash{-0.1}{-0.2}{0.7}{0.7} X_1 & \eta^{L}_{K}
  & + & \enbox{-0.1}{-0.2}{0.7}{0.7} X_2 & \tK\\
  d_2^M\rule{0cm}{0.7cm}     & \enbox{-0.1}{-0.2}{0.7}{0.7} \eta^M_N     & = &
  \enboxdash{-0.1}{-0.2}{0.7}{0.7} X_1\rule{0cm}{0.6cm} &
  \hspace{0.1cm}\enbox{-0.1}{-0.2}{0.7}{0.7} \eta^{L}_N & + &
  \enbox{-0.1}{-0.2}{0.7}{0.7} X_2 & \eta^{K}_N & + & \enbox{-0.1}{-0.2}{0.7}{0.7} X_3 & \tN.
\end{array}
\]

\bigskip
As mentioned at the end of Section~\ref{3 and higher} we will provide a way to build simple examples showing that the rest of the system is not successively solvable in the sense that the upper two equations ($\ast$), which define all possible $\eta^L_N \leftrightarrow E_1\in\ExtMod(\eta^L_K,\eta^K_N)$ cannot be solved independently from the lower two: \\
We choose $\tK = (1)\in D^{1\times 1}$ (i.e.\  $K=\coker(\tK)=0$) which implies $\eta^L_K = 0$ and $\eta^K_N = 0$ as $1$-cocycles. Then the two successive {\sc Yoneda} products $\eta^M_L \circ \eta^L_K$ and $\eta^L_K\circ\eta^K_N$ trivially vanish and therefore $\eta^M_K$ and $\eta^L_N$ exist. Since $\eta^M_K=0$ as a $1$-cocycle anyway, the situation is reducible to the case $c=2$:
\[
\xymatrix{
\mbox{$\left(\begin{array}{c|c|c|c}
    \ \tM\  & \eta^M_L\rule[-0.2cm]{0cm}{0.6cm} & 0 & \eta=? \\ \hline
     \cdot  &  \tL &  0 &  \eta^L_N\rule[-0.2cm]{0cm}{0.7cm}\\ \hline
     \cdot &  \cdot &  1 &  0 \\ \hline
     \cdot &  \cdot &  \cdot &  \tN\rule{0cm}{0.5cm}
\end{array}\right)$}
& \leadsto & 
\mbox{$\left(\begin{array}{c|c|c}
    \ \tM\ \rule{0cm}{0.5cm} & \eta^M_L & \eta=? \\ \hline
    \cdot\rule{0cm}{0.5cm} & \tL & \eta^L_N \\ \hline
    \cdot\rule{0cm}{0.5cm} & \cdot & \tN
\end{array}\right).$}
}
\]
Theorem~\ref{the thm} shows that for $\eta$ to exist the admissible choices of $\eta^L_N$ depend on the {\sc Yoneda} product with $\eta^M_L$. Choosing $\eta^L_N=0$ proves $\ExtMod(\eta^M_L,0,0)\neq\emptyset$. However, in Example~\ref{example necessary} we give two $1$-cocycles $\eta^M_L$ and $\eta^L_N$ with non-vanishing {\sc Yoneda} product. So starting with $\eta^M_L$ this choice of $\eta^L_N$ does not lead to a solution $\eta$.

\subsection{The cases $c\geq 4$ in ``coordinates''}\label{coord 4 and higher}

As for $c=3$, the cases $c \geq 4$ lead to quadratic systems. We only demonstrate this for $c=4$. Corollary~\ref{coro} applied to the relations matrix
\[
\left(\begin{array}{ccccc}
    \ \tM\  & \boxed{\eta^M_L}\rule[-0.2cm]{0cm}{0.6cm} & \eta^M_K & \eta^M_H& \eta^M_N \\
     \cdot  &  \tL & \enbox{-0.2}{-0.2}{2.63}{0.7}\eta^L_K & \eta^L_H &  \eta^L_N \rule[-0.2cm]{0cm}{0.7cm}\\
     \cdot \rule[-0.1cm]{0cm}{0.7cm} &  \cdot &  \tK & \eta^K_H &  \!\!\!\enbox{0.15}{-0.205}{0.7}{0.7}~\enbox{0.07}{-0.15}{0.6}{0.6}~\eta^K_N \\
      \cdot \rule[-0.1cm]{0cm}{0.7cm} &  \cdot & \cdot & \tH & \eta^H_N\\
     \cdot &  \cdot &  \cdot &  \cdot &  \tN\rule{0cm}{0.5cm}
\end{array}\right)
\]
leads to the following quadratic system of equations:
\kommentar{
$\begin{array}{rcl}
    d_2^{L} \left(\begin{array}{ccc} \eta^{L}_{K} &
    \eta^L_H & \eta^L_N\end{array}\right) & = & \left(\begin{array}{c|c|c}Y_1 &
    Y_2 & Y_3 \end{array}\right) \left(\begin{array}{ccc}\tK & \eta^K_H & \!\!\!\enbox{0.15}{-0.205}{0.7}{0.7}~\enbox{0.07}{-0.15}{0.6}{0.6}~\eta^K_N \\ \cdot & \tH & \eta^H_N\rule{0cm}{0.5cm}\\ \cdot & \cdot & \tN \end{array}\right)\rule[-0.8cm]{0cm}{1cm} \\
  \hline \hline 
  \setlength{\unitlength}{1cm} d_2^M \eta^M_L & = & X_1 \tL\rule{0cm}{0.7cm}\\
    d_2^{M} \left(\begin{array}{ccc} \eta^M_K & \eta^M_H & 
    \eta^M_N\end{array}\right) & = & X_1 \left(\begin{array}{ccc}\eta^L_K & \eta^L_H & \eta^L_N \end{array}\right) + \left(\begin{array}{c|c|c}
    X_2 & X_3 & X_4 \end{array}\right) \left(\begin{array}{ccc}\tK & \eta^K_H & \!\!\!\enbox{0.15}{-0.205}{0.7}{0.7}~\enbox{0.07}{-0.15}{0.6}{0.6}~\eta^K_N \\ \cdot & \tH & \eta^H_N\rule{0cm}{0.5cm}\\ \cdot & \cdot & \tN \end{array}\right).
\end{array}$
}
\[
\begin{array}{r@{\hspace{0.2cm}}lcr@{\hspace{0.2cm}}lcr@{\hspace{0.2cm}}lcr@{\hspace{0.2cm}}lcr@{\hspace{0.2cm}}l}
  d_2^K & \eta^K_H & = & \enbox{-0.15}{-0.2}{0.7}{0.7} Z_1 & \tH \\
  d_2^K\rule[-0.3cm]{0cm}{1cm} & \!\!\!\!\enbox{0.15}{-0.205}{0.7}{0.7}~\enbox{0.07}{-0.15}{0.6}{0.6}~\eta^K_N & = & \enbox{-0.15}{-0.2}{0.7}{0.7} Z_1 & \eta^H_N & + &  \enbox{-0.15}{-0.2}{0.7}{0.7} Z_2 & \tN \\
  \hline\hline
  d_2^{L}\rule{0cm}{0.7cm} & \eta^{L}_{K}    & = &
  \enbox{-0.2}{-0.2}{0.7}{0.7} Y_1 & \tK\\
  d_2^{L}\rule{0cm}{0.7cm} & \enbox{-0.1}{-0.2}{0.7}{0.7} \eta^L_H    & = &
  \enbox{-0.2}{-0.2}{0.7}{0.7} Y_1 & \eta^K_H & + & \enbox{-0.2}{-0.2}{0.7}{0.7} Y_2 & \tH\\
  d_2^{L}\rule[-0.3cm]{0cm}{1cm} & \enbox{-0.1}{-0.2}{0.7}{0.7} \eta^{L}_N & = &
  \enbox{-0.2}{-0.2}{0.7}{0.7} Y_1 & \!\!\enbox{0.15}{-0.205}{0.7}{0.7}~\enbox{0.07}{-0.15}{0.6}{0.6}~\eta^K_N & + & \enbox{-0.2}{-0.2}{0.7}{0.7} Y_2 & \eta^H_N & + &
  \enbox{-0.2}{-0.2}{0.7}{0.7} Y_3 & \tN \\
  \hline \hline
  d_2^M\rule[-0.3cm]{0cm}{0.9cm}     & \eta^M_{L}        & = &
  \encircle{0.375}{0.1}{0.7}~X_1 & \tL\\
  \hline
  d_2^M\rule{0cm}{0.7cm}     & \enbox{-0.1}{-0.2}{0.7}{0.7} \eta^M_{K} & = &
  \enboxdash{-0.1}{-0.2}{0.7}{0.7} X_1 & \eta^{L}_{K}
  & + & \enbox{-0.1}{-0.2}{0.7}{0.7} X_2 & \tK\\
  d_2^M\rule{0cm}{0.7cm}     & \enbox{-0.1}{-0.2}{0.7}{0.7} \eta^M_H     & = &
  \enboxdash{-0.1}{-0.2}{0.7}{0.7} X_1\rule{0cm}{0.6cm} &
  \hspace{0.1cm}\enbox{-0.1}{-0.2}{0.7}{0.7} \eta^L_H & + &
  \enbox{-0.1}{-0.2}{0.7}{0.7} X_2 & \eta^K_H & + & \enbox{-0.1}{-0.2}{0.7}{0.7} X_3 & \tH\\
  d_2^M\rule{0cm}{0.7cm}     & \enbox{-0.1}{-0.2}{0.7}{0.7} \eta^M_N     & = &
  \enboxdash{-0.1}{-0.2}{0.7}{0.7} X_1\rule{0cm}{0.6cm} &
  \hspace{0.1cm}\enbox{-0.1}{-0.2}{0.7}{0.7} \eta^{L}_N & + &
  \enbox{-0.1}{-0.2}{0.7}{0.7} X_2 & \!\!\enbox{0.15}{-0.205}{0.7}{0.7}~\enbox{0.07}{-0.15}{0.6}{0.6}~\eta^K_N & + & \enbox{-0.1}{-0.2}{0.7}{0.7} X_3 & \eta^H_N & + & \enbox{-0.1}{-0.2}{0.7}{0.7} X_4 & \tN.
\end{array}
\]
The middle equation $d_2^M\eta^M_L=X_1\tL$ is as always solvable and independent from the rest. An analogous argument to the one given at the end of Subsection~\ref{coord extmod3} shows that in general the remaining blocks of equations (three for $c=4$) cannot be treated independent from each other. This still leaves us with a quadratic system.

\section{Examples}\label{examples}

The following examples have been computed using $\homalg$ \cite{BR}, which was extended by the second author to include the {\sc Yoneda} equivalence and the {\sc Yoneda} product. The detailed computations and more examples can be found on the homepage of $\homalg$~\cite{homalg}. See also \cite{QR08} for explicit computations with $1$-extension modules.

\subsection{The most simple example ($c=2$)}\label{example most simple}
We illustrate the modelling theorem~\ref{action} using this simple example. Let $D=\Z$. Since $D$ is a principal ideal ring, $\Ext^2_\Z=0$ and the condition of the existence theorem~\ref{the thm} is always fullfilled. We set $M=L=N=\Z/2\Z$ and consider the associated relation matrix of a corresponding $2$-extension module
\[
  \tE := \left(\begin{array}{ccc} 2 & \eta^M_L & \eta=? \\ 0 & 2 & \eta^L_N \\ 0 & 0 & 2 \end{array}\right).
\]
Since $\Hom(\Z/2\Z,\Z/2\Z)\cong\Ext^1_\Z(\Z/2\Z,\Z/2\Z)\cong\Z/2\Z$ we conclude for the first extension group $\Ext^1(\eta^M_L,\eta^L_N)$ of $1$-cocycles $\eta^M_L,\eta^L_N\in\{(0),(1)\}$ that
\[
  \Ext^1(\eta^M_L,\eta^L_N) =\left\{ \begin{array}{rl} \Z/2\Z & ,\quad \mbox{if } 
\eta^M_L=(0) \mbox{ and } \eta^L_N = (0) \\
  0 & ,\quad \mbox{else.} \end{array} \right.
\]
This means
\[
  |\ExtMod(\eta^M_L,\eta^L_N)| =\left\{ \begin{array}{rl} 2 & ,\quad \mbox{if } \eta^M_L=(0) \mbox{ and } \eta^L_N = (0) \\
  1 & ,\quad \mbox{else,} \end{array} \right.
\]
so only the pair $(\eta^M_L,\eta^L_N)=((0),(0))$ is non-rigid.

\subsection{The condition of Theorem~\ref{the thm} is sufficient ($c=2$)}\label{example sufficient}
~ 
Let $D = \Q[x,y,z]$.
Let $\tM$, $\tL$ and $\tN$ be relation matrices for the modules
$M$, $L$ and $N$:
\[
  \tM := \left(\begin{array}{c} x \\ y \\ z \end{array}\right), \quad
  \tL = \tN := \left(\begin{array}{c} x^5 \\ z \end{array}\right).
\]
Let $0\neq\eta^M_L \in \Ext^1(M,L)$ and $0\neq\eta^L_N\in \Ext^1(L,N)$ be $1$-cocycles represented by the matrices $\tilde\eta^M_L$ and $\tilde\eta^L_N$:
\[
  \tilde\eta^M_L := \left(\begin{array}{c} 0 \\ x^4 \\ 0 \end{array}\right), \quad
  \tilde\eta^L_N := \left(\begin{array}{c} 0 \\ x \end{array}\right).
\]
A representing matrix for the \textsc{Yoneda} product of the
$1$-cocycles is
\[
  \widetilde{\eta^M_L \circ \eta^L_N} = \left(\begin{array}{c} 0 \\ x \\ 0 \end{array}\right),
\]
which is trivial in $\Ext^2(M,N)$, so we know
that $\ExtMod(\eta^M_L,\eta^L_N) \neq \emptyset$. Indeed, a particular
solution $\eta$ in (\ref{extmod2}) is
\[
 \eta = \left(\begin{array}{c} 0 \\ 1 \\ 0 \end{array}\right).
\]
and since $\Ext^1(\eta^M_L,\eta^L_N) = 0$ (even though $\Ext^1(M,N)\neq 0$) we get $|\ExtMod(\eta^M_L,\eta^L_N)|=1$ and the pair $(\eta^M_L,\eta^L_N)$ is rigid.
\subsection{The condition of Theorem~\ref{the thm} is necessary ($c=2$)}\label{example necessary}
Consider the ring $D=\Q[x,y]$. Let $\tM := \left(\begin{smallmatrix}x \\ y \end{smallmatrix}\right)$,
$\tL := \left(\begin{matrix} x & y \end{matrix}\right)$, $\tN :=
(x)$, further $\tilde\eta^M_L := \left( \begin{smallmatrix}
    1&0\\0&-1\end{smallmatrix}\right)$ and $\tilde\eta^L_N := (1)$.  In this
situation the \textsc{Yoneda} product $\eta^M_L \circ \eta^L_N \in \Ext^2(M,N)$ does
not vanish and indeed the system (\ref{extmod2}) has no solution.

\subsection{The necessary condition of Corollary~\ref{necessary} is not sufficient for $c\geq 3$}\label{example not sufficient}

Let $D=\Q[x,y]$, $\tM :=
\left(\begin{smallmatrix}x\\y \end{smallmatrix}\right)$, $\tL :=  \left(\begin{smallmatrix}x\\y \end{smallmatrix}\right)$, $\tK := (x^2 y)$ and
$\tN := (y)$. If we choose $\eta^M_L = \left( \begin{smallmatrix}
    0\\\noalign{\medskip}1\end{smallmatrix} \right)$, $\eta^L_K =
\left( \begin{smallmatrix} xy\\ 0\end{smallmatrix} \right) $
and $\eta^K_N = (x)$, both \textsc{Yoneda} products $\eta^M_L \circ
\eta^L_K$ and $\eta^L_K \circ \eta^K_N$ vanish.
However, there is no simultaneous solution for $\eta^M_K$, $\eta^L_N$
and $\eta^M_N$, so no $3$-extension module exists.

\PUSH{Appendix_Yoneda.tex}%
\input Appendix_Yoneda.tex%
\POP


%% file: Appendix_Yoneda.tex
\appendix

\section{$\Ext$'s as Satellites}\label{satellites}

Let $0 \leftarrow M \xleftarrow{d_0} P_0 \xleftarrow{d_1} P_1 \xleftarrow{d_2}
\cdots \xleftarrow{d_c} P_c \xleftarrow{d_{c+1}} P_{c+1}$ be the beginning of
a projective resolution of $M$. A $c$-cocycle is by definition a morphism
$\eta:P_c \to N$ with $d_{c+1}^*(\eta)=d_{c+1}\eta=0$. I.e.\ $\eta$ factors
over.
\begin{equation}\tag{$c$-Syz}\label{nsyz}
    P_c/\img(d_{c+1}) = P_c/\ker(d_c) \stackrel{d_c}{\cong} \ker(d_{c-1}) =:
   K_c.
\end{equation}

$K_c$ is called the $c$-th syzygies module of $M$, which is due to
{\sc Schanuel}'s Lemma uniquely defined up to \emph{projective equivalence}.
This establishes the well-known equivalence between the $c$-th derived
functor
\begin{eqnarray*}
 \RR^c\Hom(-,N)(M) &:=& \defect( \Hom_R(P_{c-1},N) \xrightarrow{d_c^*}
\Hom_R(P_c,N) \xrightarrow{d_{c+1}^*} \Hom_R(P_{c+1},N)) \\
 &=& \frac{\{\eta: P_c \to N \mid  0=d_{c+1}^*(\eta):=d_{c+1}\eta\}}{\{
d_c^*\phi:=d_c\phi \mid \phi: P_{c-1} \to N\}}
\end{eqnarray*}
and the $c$-th right \emph{satellite}\footnote{The right satellites of a
\emph{left exact} functor coincide with the right derived functors
(cf.~\cite[Theorem~V.6.1]{CE} and \cite[Prop.~IV.5.8]{HS}).}
\begin{eqnarray*}
  \Sat^c\Hom(-,N)(M) &:=& \coker( \Hom_R(P_c,N)
\xrightarrow{d_c^*} \Hom_R(K_c,N) ) \\
  &=& \{\eta:K_c \to N \} / \{ d_c^*\phi=d_c\phi \mid \phi: P_{c-1} \to
N\}
\end{eqnarray*}
(cf.~\cite[Section III.2, Prop.~IV.5.8, Exercises IV.7.3 and
IX.3.1]{HS} and \cite[III.(6a),(6a')]{CE}). In words, the first right
satellite of $\Hom_R(-,N)$ applied to $M$ is the abelian group of
all morphisms $K_c\to N$, modulo those which factor over $P_{c-1}$ (i.e.\
 which extend to $P_{c-1}$).

\section{The {\sc Yoneda} Composite of Extensions}\label{composition}

Two $1$-extensions or just \emph{extensions} $G$ and $G'$ of $N$ with $M$
are called \emph{congruent} $G\equiv G'$ if there exists a chain map $(\id_M,
\beta_0, \id_N) :
G \rightarrow G'$, i.e. a commutative diagram
\begin{equation}\tag{Cong}\label{cong}
\xymatrix{G: & 0 & \ar[l] M \ar@{=}[d] & G_0 \ar[l]
\ar^{\beta_0}[d] & N \ar@{=}[d] \ar[l] & 0 \ar[l] \\
  G': & 0 & \ar[l] M & \ar[l] G_0' & \ar[l] N & \ar[l]  0.
}
\end{equation}
In this case, the Five Lemma shows that the middle homomorphism $\beta_0$ is
an isomorphism, hence congruence of extensions is a reflexive, symmetric,
and transitive relation.

For a short exact sequence 
$$\xymatrix{
  G : 0 & \ar[l] M & \ar[l] G_0 & \ar[l]_{\iota} N & \ar[l]  0
}$$
and a morphism $\phi : N \rightarrow N'$ the following diagram with the
\emph{pushout} $G'_0$ of $G_0\xleftarrow{\iota} N \xrightarrow{\phi} N'$
(cf.~\cite[Exercise II.9.2]{HS}) is commutative:
\begin{equation}\tag{Pushout}\label{push}
\xymatrix{
  G : 0 & \ar[l] M \ar@{=}[d] & \ar[l] G_0 \ar[d] & \ar_{\iota}[l] N
  \ar^{\phi}[d] & \ar[l]  0 \\
  G' : 0 & \ar[l] M & \ar[l] G_0' & \ar[l] N' & \ar[l]  0.
}
\end{equation}
Recall, the pushout is the cokernel of
\[
  N \xrightarrow{\left(\begin{array}{cc} \iota & \phi
\end{array}\right)} G_0 \oplus N' \rightarrow G_0' \to 0.
\]
The resulting exact sequence $G'$ is called the \emph{{\sc Yoneda} composite} of $G$ and $\phi$ and denoted by $G \phi$, it is unique up to congruence.

Likewise, the \emph{{\sc Yoneda} composite} $G'' = \vartheta G$ with a morphism $\vartheta : M'
\rightarrow M$ is the pullback $G'_0$ of $M' \xrightarrow{\vartheta} M
\xleftarrow{\pi} G_0$:
\begin{equation}\tag{Pullback}\label{pull}
\xymatrix{
  G : 0 & \ar[l] M & \ar_{\pi}[l] G_0 & \ar[l] N & \ar[l]  0 \\
  G'' : 0 & \ar[l] M' \ar^{\vartheta}[u]& \ar[l] G'_0 \ar[u] & \ar[l] N \ar@{=}[u] & \ar[l]  0.
}
\end{equation}
Recall, the pullback is the kernel of
\[
  0\to G_0' \to G_0 \oplus M' \xrightarrow{\left(\begin{array}{c} \pi \\ -\vartheta
\end{array}\right)} M.
\]

Write a $c$-extension $G$ as the composite of $c$ short exact sequences (cf.~Section~\ref{expository}):
\[
  G = G_0 \circ G_1 \circ \ldots \circ G_{c-1}.
\]
The \emph{{\sc Yoneda} composite} of this $c$-extension with morphisms is defined as
\begin{equation}\tag{$G\phi$}\label{Gphi}
  G \phi := G_0 \circ G_1 \circ \ldots \circ (G_{c-1}\phi)
\end{equation}
and
\begin{equation}\tag{$\vartheta G$}
  \vartheta G := (\vartheta G_0) \circ G_1 \circ \ldots \circ G_{c-1}.
\end{equation}

For $c > 1$, congruence is a wider relation. Write any $c$-extension $G$ as the composite of $c$ short exact sequences:
$$ G = G_0 \circ G_1 \circ \ldots \circ G_{c-1}.$$
Then congruence of $c$-extensions is defined as follows (cf.~\cite[Section~III.5]{ML}):
\begin{itemize}
\item $G_0 \circ G_1 \circ \ldots \circ G_{c-1} \equiv G'_0 \circ G'_1 \circ \ldots \circ G'_{c-1}$ if $G_i \equiv G_i'$ for all $i$;
\item $G_0 \circ \ldots \circ (G_{i-1} \beta) \circ G_i \circ \ldots \circ G_{c-1}
\equiv G_0 \circ \ldots \circ G_{i-1} \circ (\beta G_i) \circ \ldots \circ
G_{c-1}$ for a matching homomorphism $\beta$.
\end{itemize}

This congruence relation defines a module for every $c \geq 1$:
\[ 
  \YExt^c(M,N) := \{\text{$c$-extensions of $N$ with $M$} \}/\equiv.
\]
For $c=0$ one sets $\YExt^0(M,N) := \Hom(M,N)$. The natural equivalence of functors
\[
  \xymatrix@1{\Yon^c:\Ext^c(-,-) \ar[r] & \YExt^c(-,-)}
\]
is explicitly described in Appendix~\ref{equivalence}.

\kommentar{
The {\sc Yoneda} equivalence allows avoiding both the use of injective resolutions  in the construction of $\Ext^c(M,-)$ as the right derived functor $\RR^c \Hom(M,-)$ and the use of projective resolutions for $\Ext^c(-,N)=\RR^c \Hom(-,N)$. This establishes the existence of such functors independent of the existence of enough injectives or projectives in the relevant category.}

\section{From $c$-Extensions to $c$-Cocycles and back}\label{equivalence}

To a $c$-extension or shortly an extension of $N$ by $M$ corresponds a cocycle in
$\Ext^c(M,N)$. To describe this we need the following two ingredients, which both appeared in Appendix \ref{satellites}. First, the representation of the $c$-cocycle $\eta$ by a map $K_c\to N$ (recall, $K_c$ is defined up to projective equivalence). Second, the $c$-extension
\[
  G^c_M: 0\leftarrow M \leftarrow P_0 \leftarrow \cdots \leftarrow P_{c-1} \leftarrow K_c \leftarrow 0,
\]
depending only on $M$ up to projective equivalence. The natural equivalence $\Yon^c$ satisfies the defining equation
\[
  G^c_M \  \eta = \Yon^c(\eta),
\]
where $G^c_M \  \eta$ is the {\sc Yoneda} composite defined above (see (\ref{Gphi})). The inverse of $\Yon^c$ is given by lifting the identity map $\id_M$ to $\eta:K_c\to N$:
\begin{equation}\tag{$\mathrm{Lift}^c$}
  \xymatrix{
    0 & \ar[l] M \ar@{=}[d] & \ar[l] P_0 \ar@{.>}[d]^{\eta_0} & \cdots \ar[l] & P_{c-1} \ar[l] \ar@{.>}[d]^{\eta_{c-1}} &\ar[l]_{d_c}
  K_c \ar@{.>}[d]^{\eta} & \ar[l]   0 \\
    0 & \ar[l] M & \ar[l] G_0 & \cdots \ar[l] & G_{c-1} \ar[l] & \ar[l] N & \ar[l] 0.
  }
\end{equation}

See \cite[Theorem III.2.4]{HS} for $c=1$ and  \cite[Theorem IV.9.1]{HS} and \cite[Theorem III.6.4]{ML} for the general case.

\section{$c$-Cocycles and the \textsc{Yoneda} Product}\label{lift}

Let $\eta^M_L \in \Ext^c(M,L)$ and $\eta^L_N \in \Ext^{c'}(L,N)$ be cocycles for
$c,c' \geq 0$, $0 \leftarrow M \leftarrow P$ and $0 \leftarrow L \leftarrow Q$ projective resolutions of
$M$ and $L$.
Lift $\eta^M_L$ to $\eta$ as in the diagram
\[
\xymatrix{
  & & P_c \ar[dl]_{\eta^M_L} \ar@{.>}^{\eta_0}[d] & P_{c+1} \ar[l] \ar@{.>}^{\eta_1}[d] & \ldots \ar[l] & P_{c+c'-1} \ar[l] \ar@{.>}^{\eta_{c'-1}}[d] & P_{c+c'} \ar[l] \ar@{.>}_{\eta}[d] \ar@{-->}@/^1.5pc/^{\eta^M_L \circ \eta^L_N}[dd]\\
  0 & L \ar[l] & Q_0 \ar[l] & Q_1 \ar[l] & \ldots \ar[l] & Q_{c'-1} \ar[l] &
  Q_{c'} \ar[l] \ar_{\eta^N_L}[d]\\
  & & & & & & \tL
}
\]
then the composite of morphisms $\eta\,\eta^N_L$ is an element of
$\Ext^{c+c'}(M,N)$. This defines a product
\[
\begin{array}{rclcc}
  \Ext^c(M,L) & \times & \Ext^{c'}(L,N) & \to & \Ext^{c+c'}(M,N) \\
  (\eta^M_L &,& \eta^L_N) & \mapsto & \eta^M_L \circ \eta^L_N := \eta\,
  \eta^N_L
\end{array}
\]
which is called the \emph{{\sc Yoneda} product} of cocycles.

\section{Triangular systems over rings}\label{appendix triangular}

It is a well known fact that a system in triangular shape over a \emph{field} can be solved successively. We call systems with this feature ``triangular systems''. Over rings triangular-shape systems need not in general be triangular. For example let $D=\Z$ and consider the system
\[
\begin{array}{ccccccc}
  2x & & + & z & = & 0 \\
     & y & + & z & = & 1.
\end{array}
\]
This system is in a triangular \emph{shape}, but it cannot be solved successively. For example $y=0,z=1$ is a valid solution for the bottom equation, but for this choice the top equation is not solvable. The system is solvable, however, as we can see by setting $y=1, z=0$, and $x=0$. Hence, this set of equations over the ring $D=\Z$ is not triangular.


%% file: yoneda.bbl
\def\cprime{$'$} \def\cprime{$'$}
\providecommand{\bysame}{\leavevmode\hbox to3em{\hrulefill}\thinspace}
\providecommand{\MR}{\relax\ifhmode\unskip\space\fi MR }
\providecommand{\MRhref}[2]{%
  \href{http://www.ams.org/mathscinet-getitem?mr=#1}{#2}
}
\providecommand{\href}[2]{#2}